\documentclass[a4paper,11pt]{amsart}
\usepackage[english]{babel}
\usepackage[latin1]{inputenc}
\usepackage[T1]{fontenc}
\usepackage{amsmath}
\usepackage{amsfonts}
\usepackage{amssymb}
\usepackage{amsthm}
\usepackage{mathrsfs}
\usepackage{enumitem}
\usepackage{graphicx}
\usepackage{footnote}
\usepackage{hyperref,color}

\renewcommand{\geq}{\geqslant}
\renewcommand{\leq}{\leqslant}
\newtheorem{thm}{Theorem}

\newtheorem{rem}[thm]{Remark}
\newtheorem{cor}[thm]{Corollary}
\newtheorem{prop}[thm]{Proposition}
\newtheorem{lem}[thm]{Lemma}
\usepackage{color}
\definecolor{darkgreen}{rgb}{0,0.4,0}
\definecolor{MyDarkBlue}{rgb}{0,0.08,0.50}
\definecolor{BrickRed}{rgb}{0.65,0.08,0}

\hypersetup{
colorlinks=true,       
    linkcolor=blue,          
    citecolor=red,        
    filecolor=BrickRed,      
    urlcolor=darkgreen        
}

\title[Passage time from four to two blocks of opinions in~the voter model]{Passage time of a random walk in the quarter plane for opinions in the voter model}

\author{Irina Kurkova\and Kilian Raschel}

 \thanks{I.\ Kurkova: Laboratoire de Probabilit\'es et
        Mod\`eles Al\'eatoires, Universit\'e Pierre et Marie Curie,
        4 Place Jussieu, 75252 Paris Cedex 05, France. Email: \url{Irina.Kourkova@upmc.fr}}

 \thanks{K.\ Raschel: CNRS and Laboratoire de Math\'ematiques et Physique Th\'eorique,
Universit\'e de Tours,
Parc de Grandmont, 37200 Tours, France.
       Email: \url{Kilian.Raschel@lmpt.univ-tours.fr}
     }

 \date{\today}

\begin{document}

\begin{abstract}
A random walk in ${\bf Z}_+^2$ spatially homogeneous in the interior, absorbed at the axes, starting from an arbitrary point $(i_0,j_0)$ and with step probabilities drawn on Figure \ref{v-Figfig} is considered. The trivariate generating function of probabilities that the random walk hits a given point $(i,j) \in {\bf Z}_+^2 $ at a given time $k\geq 0$ is made explicit. Probabilities of absorption at a given time $k$ and at a given axis are found, and their precise asymptotic is derived as the time $k\to \infty$. The equivalence of two typical ways of conditioning this random walk to never reach the axes is established. The results are also applied to the analysis of the voter model with two candidates and initially, in the population ${\bf Z}$, four connected blocks of same opinions. Then, a citizen changes his mind at a rate proportional to the number of its neighbors that disagree with him. Namely, the passage from four to two blocks of opinions is studied.

\medskip

\medskip


\medskip

\medskip

\end{abstract}

\maketitle

\section{Introduction}
\label{v-Intro}

\subsection*{Context}

   Random walks with small steps in the quarter plane ${\bf Z}_+^2=\{0,1,2,\ldots \}^2$ spatially homogeneous in the
   interior and on each of the two axes are now rather well studied.
   The analytic approach \cite{FIM} elaborated by Fayolle, Iasnogorodski and Malyshev
   provided the generating function, say $H(x,y)$,
    of the stationary probabilities
    in the ergodic case, and also
   that of the Green functions
   in the transient case.
   Further analysis allowed to
   compute the asymptotic of these quantities along any path in ${\bf Z}_+^2$, see \cite{KM,KR,MALY,KRSp4,thesis}.


{The main motivation of the present work is to
develop  a method for incorporating the parameter $z$ of time into
this approach, in order to derive the trivariate generating function
$H(x,y;z)$ of the probabilities $h_{i,j;k}$ that the walk is in
state $(i,j)$ at time $k$, and to obtain asymptotic results from it.
Being able to deal with this additional time variable $z$ is
actually important in combinatorics (e.g., to count certain numbers
of walks confined to the quarter plane, see \cite{Ra}) and in
probability as well (e.g., to compute the distribution of some
hitting times). This is one of the few attempts is that direction,
after \cite{Blanc,FIMit}. The second motivation of this article is
the application that it has for the voter model: indeed, it
completes results of \cite{Men01,Men08} about the hitting time of
the so-called Heaviside configuration in the voter model with
initially four blocks of opinions.}

\subsection*{Voter model}

By the voter model we mean a continuous-time process on $\{0,1\}^{{\bf Z}}$ (here and throughout, ${\bf Z} = \{\ldots ,-1,0,1,\ldots \}$) that can
be interpreted as follows: initially, at each site of ${\bf Z}$, there is zero or one particle;
then a particle appears (resp.\ disappears) at an empty (resp.\ occupied) site $x$ according to
an exponential law with a rate proportional to the number of nearest neighbors of $x$ which are occupied
(resp.\ empty). Moreover, we assume that the initial state appertains to the set of configurations
having a finite number of empty (resp.\ occupied) sites on the left (resp.\ right) of the origin $0$, 
{see \eqref{v-config} for an example.}
In particular, this implies that at any time the process will belong to this set of configurations.
As a consequence there is, at any time, a finite number of ``$01$'' (resp.\ ``$10$''), i.e., a
finite number of pairs of sites $(x,x+1)$ with zero (resp.\ one) particle at $x$ and one
(resp.\ zero) particle at $x+1$.

The underlying discrete-time voter model is a Markov chain with the following
dynamic: denote by $\mathscr{C}_{k}$ the configuration at time $k$ (and remember
that according to the previous paragraph, there is only a finite number of ``$01$''
and ``$10$'' in $\mathscr{C}_{k}$); next, in order to construct $\mathscr{C}_{k+1}$,
one first chooses with a uniform distribution one of these ``$01$'' and ``$10$'' in
$\mathscr{C}_{k}$, then one replaces it, with probability $1/2$, by ``$00$'' or ``$11$''.

If the voter model starts from the Heaviside configuration, i.e., the configuration
having only occupied (resp.\ empty) sites on the left (resp.\ right) of the origin, then at
any time, the process will be a translation of it. This fact suggests to consider the following
equivalence relation: two configurations are equivalent if they are translations
the one of the other.

From now on, we shall work on the underlying quotient space, the equivalence classes of
which being identified by finite sets of positive integers $(X_{1},Y_{1},\ldots , X_{N},Y_{N})$:
     \begin{equation}
     \label{v-config}
          \ldots 111
                 \overbrace{000000}^{X_{1}}
                 \overbrace{11111}^{Y_{1}}
                 \overbrace{00000}^{X_{2}}
                 \overbrace{11111}^{Y_{2}}
          \ldots \ldots
                 \overbrace{000}^{X_{N}}
                 \overbrace{1111}^{Y_{N}}
                 000
          \ldots ,
     \end{equation}
$N$ being the number of finite blocks of zeros (or ones) and $X_{\ell}$ (resp.\ $Y_{\ell}$),
$\ell\in\{1,\ldots ,N\}$, the size of the $\ell$th block of zeros (resp.\ ones). The number
$N$ of finite blocks of zeros is a non-increasing function of the time;
furthermore, $N=0$ corresponds to the class of the Heaviside configuration.

We refer to \cite{Lig} for additional details on the voter model
and, more generally, for further information about interacting
particle systems.

\subsection*{Hitting time of the Heaviside configuration}

Let $\tau$ denote the hitting time of the Heaviside configuration. It is proved in \cite{Men01} that for any initial configuration,
     \begin{align}
     \label{ttt1}
          {\bf E} [\tau^{3/2-\epsilon}]<\infty, \quad \forall \epsilon>0,
          \\
          \label{ttt2}{\bf E}[\tau^{3/2+\epsilon}]=\infty,\quad \forall \epsilon>0.
     \end{align}
Statement \eqref{ttt1} is proved by an adequate use of Lyapunov functions. To show \eqref{ttt2}, it suffices to do it only for initial states with $N=1$ in \eqref{v-config}; that is done in \cite{Men01}, by applying results on passage time moments proved in \cite{ASP}.

 With the notations \eqref{v-config}, consider the process
  starting from  an initial state with $N=1$:
$(X_{1},Y_{1})=(X_{1}(k),Y_{1}(k) )_{k\in{\bf Z}_+}$. We rename it here
$(X,Y)=(X(k),Y(k) )_{k\in{\bf Z}_+}$; we have:
     \begin{equation*}
          \ldots 111
                 \overbrace{000000}^{X}
                 \overbrace{11111}^{Y}
                 000
          \ldots .
     \end{equation*}
The process $(X,Y)$ is a Markov chain on ${\bf Z}_{+}^{2}$ which is absorbed
as it reaches the boundary, since the Heaviside configuration is an
absorbing state for the voter model. Moreover, using the dynamic of
the discrete-time voter model explained above, we notice that
$(X,Y)$ has homogeneous transition probabilities in the interior of
${\bf Z}_{+}^{2}$ equal to (with obvious notations)
\begin{equation*}
p_{1,0}=p_{1,-1}=p_{0,-1}= p_{-1,0}=p_{-1,1}=p_{0,1}=1/6
\end{equation*}
 and the
others to $0$, see Figure \ref{v-Figfig}. Further, the hitting time $\tau$ can be expressed as
     \begin{equation}
     \label{def_ht}
          \tau=\inf\{k\in{\bf Z}_+: X(k)=0 \text{ or } Y(k)=0\}.
     \end{equation}
Define also the hitting times of the horizontal and vertical axes:
     \begin{equation}\label{v-def_hitting_times}
          S=\inf\{k\in{\bf Z}_+: Y(k)=0\},\quad
          T=\inf\{k\in{\bf Z}_+: X(k)=0\},
     \end{equation}
so that $\tau=\inf\{S,T\}$.


\subsection*{Main results}

{The present work is constituted by four main points, that we now describe.}

The first one is a direct consequence of \cite{thesis} and concerns explicit expressions for the probabilities that the process is absorbed at some site of the boundary in a given time, namely, ${\bf P}_{(i_{0},j_{0})}[S=k]$ and ${\bf P}_{(i_{0},j_{0})}[T=k]$, for any $k\in{\bf Z}_+$ and any
$(i_{0},j_{0})\in {\bf Z}_+^2$. For this we shall use Proposition \ref{explicit_integral_marches_en} of Section
\ref{DistributionTheVoterModel}, taken from \cite[Chapter F]{thesis}, which gives an integral representation of the generating
functions
     \begin{eqnarray}
               h^{i_0,j_0}(x;z)  &=&
               \sum_{i\geq 1}\sum_{k\geq 0}
               {\bf P}_{(i_{0} , j_{0})}
               [(X,Y)\text{ hits } (i,0)\text{ at time $k$}]
               x^{i} z^{k},\label{v-def_generating_functions_x}
               \\
               \widetilde{h}^{i_0,j_0}(y;z)  &=&
               \sum_{j \geq 1}\sum_{k\geq 0} {\bf P}_{(i_{0} , j_{0})}
               [(X,Y)\text{ hits }(0,j)\text{ at time $k$}]
                y^{j} z^{k}\label{v-def_generating_functions_y}.
     \end{eqnarray}
  Then ${\bf P}_{(i_{0},j_{0})}[S=k]$ and
${\bf P}_{(i_{0},j_{0})}[T=k]$ can be expressed from $h^{i_0,j_0}(1;z)$ and
$\widetilde h^{i_0,j_0}(1;z)$ via the Cauchy formul\ae. Note, besides, that we also find the trivariate function
  \begin{equation}
  \label{gf}
  H^{i_0,j_0}(x,y;z) = \sum_{i,j\geq 1}\sum_{k\geq 0}
          {\bf P}_{(i_{0},j_{0})}
          [ (X(k),Y(k))=(i,j)] x^{i-1} y^{j-1} z^{k}
  \end{equation}
thanks to the functional
equation \eqref{v-functional_equation}.
{For the voter model,
this means that we find explicit expressions for the probabilities that the process hits the Heaviside configuration at any fixed time, with the additional information of the size of the blocks at the time of absorption.}

\begin{figure}[t]
\begin{center}
\begin{picture}(000.00,780.00)
\hspace{-100mm}\includegraphics{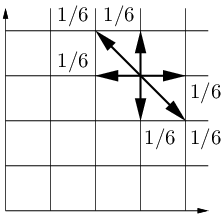}
\end{picture}
\end{center}
\vspace{-240mm}
\caption{Transitions of the process $(X,Y)$ {in the interior of the quarter} {plane ${\bf Z}_+^2$; on the boundary, the process is absorbed}}
\label{v-Figfig}
\end{figure}

{The second point is our main result, and is new,
 to the best of our knowledge. It is about the asymptotic tail distribution of the hitting time $S$.}

\begin{thm}
\label{main_result_group_six_drift_zero}
As the time $k\to \infty$, we have
     \begin{equation}
     \label{zghi}
          {\bf P}_{(i_{0},j_{0})}[S=k]
          = \frac{9}{16}\left(\frac{3}{\pi}\right)^{1/2}
          \frac{i_{0} j_{0}(i_{0}+j_{0})}{k^{5/2}}(1+o(1)).
     \end{equation}
     \end{thm}
{Theorem \ref{main_result_group_six_drift_zero}
will be a  consequence of Proposition \ref{v-lemma_h1} and of
classical singularity analysis \cite[Sections 6.2--6.4]{FLAJ}, see
Section \ref{AsymptoticTheVoterModel}. Theoretically, the methods
developed in this paper could work for (the hitting time of the
boundary of) any random walk with transitions to the eight nearest
neighbors, see Remark \ref{gluing}.}

We now introduce the third  point of our work. First notice that the
transition probabilities of the walk are such that ${\bf
P}_{(i_{0},j_{0})}[S=k]={\bf P}_{(j_{0},i_{0})}[T=k]$, see Figure
\ref{v-Figfig}, and that the quantity $i_{0} j_{0}(i_{0}+j_{0})$ in
\eqref{zghi} is invariant by $(i_0,j_0)\mapsto (j_0,i_0)$.
Accordingly, the asymptotic of ${\bf P}_{(i_{0},j_{0})}[T=k]$ is
exactly the same as that of ${\bf P}_{(i_{0},j_{0})}[S=k]$. Further,
\begin{equation}
\label{sum}
     {\bf P}_{(i_{0},j_{0})}[\tau=k]={\bf P}_{(i_{0},j_{0})}[S=k]+{\bf P}_{(i_{0},j_{0})}[T=k],
\end{equation}
so that Theorem \ref{main_result_group_six_drift_zero} entails
the following corollary.
\begin{cor}
\label{cor2}
The result \eqref{ttt2} proved in \cite{Men01} for $\epsilon>0$ also holds for $\epsilon=0$.
\end{cor}
{The latter completes the results of \cite{Men01,Men08}. It is worth noting that this corollary can also be obtained as a consequence of results in \cite{VA1,DW}. More details are provided at the end of this introduction.}

Finally, the last result in our paper is the following: the precision of the asymptotic result \eqref{zghi} implies
Corollary~\ref{dodo} below that compares two typical ways on
conditioning the process $(X,Y)$ to never reach the axes. Define
$h(i_{0},j_{0})=i_{0}j_{0}(i_{0}+j_{0})$ (in fact, using methods closed to \cite{KRSp4}, it could be proved that it is the unique positive
harmonic function associated with the walk $(X,Y)$ absorbed at the
boundary of ${\bf Z}_+^2$).
\begin{cor}
\label{dodo} The Doob $h$-process of $(X,Y)$ coincides in
distribution with the limit, as $k\to \infty$, of the process
conditioned on $\{\tau>k\}$.
\end{cor}

The proof relies on the following precise asymptotic, as $k\to
\infty$,
     \begin{equation*}
          {\bf P}_{(i_{0},j_{0})}[\tau \geq k ]
          = \frac{27}{16} \left(\frac{3}{\pi}\right)^{1/2} \frac{i_{0} j_{0}(i_{0}+j_{0})}{k^{3/2}}(1+o(1)),
     \end{equation*}
which is a direct consequence of Theorem
\ref{main_result_group_six_drift_zero} and \eqref{sum}. It is then
carried out  by a standard reasoning as in \cite{KRSp4} or \cite[Chapter
F]{thesis}.

\subsection*{Other approaches}
We close this introduction by mentioning other possible approaches
for analyzing asymptotic tail distribution of hitting times for
random walks in cones of ${\bf Z}^d$. First, as already quoted,
methods using Lyapunov functions in \cite{ASP,Men01,Men08} show the
finiteness or infiniteness of hitting times' moments. {A series of
tail distribution estimates for
 hitting times is presented in \cite{VA1}, by using potential theory. They result in upper and lower bounds for ${\bf P}_{(i_{0},j_{0})}[\tau \geq k ]$, which are enough to obtain Corollary \ref{cor2} in this paper.
  In a  recent work \cite{DW},
   the tail asymptotic of the hitting time up to a multiplicative
factor is obtained by comparison with Brownian motion. This is another way to deduce Corollary \ref{cor2}.}
  All these methods are powerful for rather general random walks in conic
  domains of ${\bf Z}^d$, but do not give as much accurate results
 as Theorem \ref{main_result_group_six_drift_zero} in this paper.
       Finally, let us mention the paper \cite{KP}, where
an approach based on an extension of the Karlin-McGregor formula is
applied to the family of the so-called non-colliding random walks.
The latter leads to precise results, but it exploits a particular
independence property of this family, and therefore seems to be
restricted to this class of models.

\section{Exact distribution of the hitting times of both axes}
\label{DistributionTheVoterModel}
\setcounter{equation}{0}

{This section contains preliminary material, which is needed for Section \ref{AsymptoticTheVoterModel}, where we prove our main results.}

\subsection*{A functional equation and the kernel of the walk}
With the notations \eqref{v-def_generating_functions_x}, \eqref{v-def_generating_functions_y} and \eqref{gf} of Section \ref{v-Intro}, we can state
on $\{(x,y,z)\in {\bf C}^{3}: |x|,|y|,|z|\leq 1 \}$ (here and throughout, ${\bf C}$ denotes the complex plane)
the following crucial functional equation:
     \begin{equation}\label{v-functional_equation}
          K(x,y;z)H^{i_0,j_0}(x,y;z) =
          h^{i_0,j_0}(x;z)+\widetilde{h}^{i_0,j_0}( y;z) -x^{i_{0}} y^{j_{0}},
     \end{equation}
where $K(x,y;z)$ is the following polynomial---called the kernel of the walk---,
depending only on the walk's transition probabilities:
     \begin{equation}\label{v-def_Q}
          K(x,y;z)= x y z\textstyle[\sum_{-1\leq i,j\leq 1 }
          p_{i,j }x^{i} y^{j} -1/z].
     \end{equation}
For $z=0$, Equation \eqref{v-functional_equation} simply
becomes ${\bf P}_{(i_{0},j_{0})}[ (X(0),Y(0))=(i_{0},j_{0}) ]=1$.
For $z=1$, it becomes a functional equation between the Green functions
generating function and the absorption probabilities
generating functions; the latter is studied in \cite{KR,KRSp4,thesis}.
For the proof of \eqref{v-functional_equation}, we exactly use
the same arguments as in \cite[Chapter F]{thesis}.

We now study the set of the zeros of the kernel $K(x,y;z)$ defined in \eqref{v-def_Q}.
For this we start by remarking that it can be written alternatively
     \begin{equation}\label{v-def_q_alternative}
          K(x,y;z) = a(x;z) y^{2}+ b(x;z) y
          + c(x;z) = \widetilde{a}(y;z) x^{2}+
          \widetilde{b}(y;z) x + \widetilde{c}(y;z),
     \end{equation}
where
\begin{equation*}
     a(x;z) = z(x+1)/6, \quad b(x;z)= z x^{2}/6-x+z/6,\quad c(x;z) = z x (x+1)/6,
\end{equation*}
and
\begin{equation*}
     \widetilde{a}(y;z) = z(y+1)/6,\quad \widetilde{b}(y;z)= z y^{2}/6-y+z/6,\quad \widetilde{c}(y;z) = z y (y+1)/6.
\end{equation*}
Next, we introduce the algebraic function $Y(x;z)$ defined by
$K(x,Y(x;z);z)=0$. Note that $K(x,y;z)=0$ is equivalent to
$[b(x;z)+2a(x;z)y]^{2}=d(x;z)$, where
\begin{equation*}
d(x;z)=b(x;z)^{2}-4a(x;z)c(x;z),
\end{equation*}
so that the construction of the function
$Y(x;z)$ is equivalent to that of the square roots of the polynomial $d(x;z)$, namely, $\pm d(x;z)^{1/2}$.
For this we need the following:
     \begin{lem}\label{v-lemma_branched_points}
          Let $z\in ]0,1[$. The four roots of $x\mapsto d(x;z)$
          are
          positive and mutually distinct.
          We call them
          $x_{1}(z)<x_{2}(z)<x_{3}(z)<x_{4}(z)$.
          They satisfy $x_{1}(z)x_{4}(z)=
          x_{2}(z)x_{3}(z)=1$.
          In particular, $x_{1}(z),x_{2}(z)\in]0,1[$ and
          $x_{3}(z),x_{4}(z)\in]1,\infty [$. Further, $x_{1}(0)=x_{2}(0)=0$,
          $x_{3}(0)=x_{4}(0)=\infty $, $x_{2}(1)=x_{3}(1)=1$
          and $x_{1}(1)=7-4\sqrt{3}$, $x_{4}(1)=7+4\sqrt{3}$.
     \end{lem}

\begin{proof}
As we can easily verify, the polynomial $d(x;z)$ is reciprocal, in
other words it satisfies $x^{4}d(1/x;z)=d(x;z)$. This property allows us to write
it as a second degree polynomial
in the variable $x+1/x$. Following this way we obtain the explicit expression
of its roots: if $s_{1}(z)=3/z+1$ and $s_{2}(z)=(6/z+3)^{1/2}$, then
$x_{1}(z)=s_{1}(z)+s_{2}(z)+[(s_{1}(z)+s_{2}(z))^{2}-1]^{1/2}$ and
$x_{2}(z)=s_{1}(z)-s_{2}(z)+[(s_{1}(z)-s_{2}(z))^{2}-1]^{1/2}$,
$x_{3}(z)=1/x_{2}(z)$ and $x_{4}(z)=1/x_{1}(z)$.
All the properties of Lemma \ref{v-lemma_branched_points} immediately
follow from these explicit expressions.
\end{proof}
There are two branches of the square root of $d(x;z)$.
Each determination leads to a single-valued and meromorphic function on the
complex plane ${\bf C}$ appropriately cut, that is, in our case,
on ${\bf C}\setminus ([x_{1}(z),x_{2}(z)] \cup [x_{3}(z),x_{4}(z)])$.
We have
     \begin{equation*}
     \label{v-explicit_expression_Y}
          Y(x;z) = \frac{ -b(x;z)\pm
          d(x;z)^{1/2}}{2 a(x;z) },
     \end{equation*}
and we fix the notations of the branches $Y_0(x;z)$ and $Y_1(x;z)$ by (arbitrarily) choosing that $|Y_0(x;z)|<|Y_1(x;z)|$ on the whole of ${\bf C}\setminus ([x_{1}(z),x_{2}(z)] \cup [x_{3}(z),x_{4}(z)])$.
For more details about the
construction of algebraic functions, see, e.g., \cite{SG2}.

In a similar way, the functional
equation \eqref{v-functional_equation} defines an algebraic
function $X(y;z)$.  But it turns out that
$K(x,y;z)=K(y,x;z)$, see \eqref{v-def_Q}, so that $X(y;z)=Y(y;z)$;
in particular, all properties proved for $Y(x;z)$
 immediately result in similar ones for $X(y;z)$.

\subsection*{Explicit expression of distributions} They are obtained in the result that follows.

\begin{prop}
\label{explicit_integral_marches_en}
The function $h^{i_{0},j_{0}}(x;z)$ is equal to
     \begin{eqnarray}
          h^{i_{0},j_{0}}(x;z)\hspace{-1mm}&=&\hspace{-1mm}
          x^{i_{0}}Y_{0}(x;z)^{j_{0}}\nonumber\\\hspace{-1mm}&+&\hspace{-1mm}
          \int_{x_{1}(z)}^{x_{2}(z)}t^{i_{0}}\mu_{j_{0}}(t;z)
          \left[\frac{\partial_{t}w(t;z)}{w(t;z)-w(x;z)}-
          \frac{\partial_{t}w(t;z)}{w(t;z)-w(0;z)}\right]
          [-d(t;z)]^{1/2}\textnormal{d}t,\label{generic_equality}
     \end{eqnarray}
where
     \begin{equation}
     \label{v-def_mu}
          \mu_{j_{0}}(t;z) =  \frac{1}{[2a(t;z)]^{j_{0}}}
          \sum_{k=0}^{(j_{0}-1)/2} \binom{2k+1}{j_0}d(t;z)^{k} [
          - b(t;z)]^{j_{0}-(2k+1)},
     \end{equation}
and
     \begin{equation}
     \label{v-eq_explicit_form_CGF}
               w(t;z)=\frac{t(1+t)}
               {(t-x_{2}(z))(t-x_{3}(z)^{1/2})^{2}}.
     \end{equation}
\end{prop}

\begin{rem}
\label{gluing}
{\rm
{Equations \eqref{generic_equality} and \eqref{v-def_mu} are obtained in \cite[Chapter F]{thesis}, while \eqref{v-eq_explicit_form_CGF} is found in \cite{Ra}. It is worth noting that Equations \eqref{generic_equality} and \eqref{v-def_mu}
are valid not only for the random walk under consideration in this paper, but for all random walks with jumps to the eight nearest neighbors, see \cite{FIM} and \cite[Chapter F]{thesis}. }

{On the other hand, finding an expression for the function $w(t;z)$ happens to be quite complex in general, and dependent on the particular model. Further, in general, there is no reason for this function to be rational (in $t$) as in \eqref{v-eq_explicit_form_CGF}, or even algebraic. To be complete, we note that for our model, $w(t;z)$ is rational because a certain group of automorphisms is finite. We refer to \cite{FI,FIM,KR,Ra,KRSp4,thesis} for any details on this group, and more generally on how finding expressions as in Proposition \ref{explicit_integral_marches_en}.}

{The algebraicity of $w(t;z)$ happens to be crucial for the proof of Theorem \ref{main_result_group_six_drift_zero} (see Section \ref{AsymptoticTheVoterModel}), and this is why this article focuses on one particular model. It remains an open problem to determine, for all random walks with jumps to the eight nearest neighbors, the asymptotic tail distribution of the hitting time of the boundary, by using analytic methods as in this paper.}}
\end{rem}

Now, using the partial fraction expansion (direct consequence of \eqref{v-eq_explicit_form_CGF})
     \begin{multline*}
          \frac{\partial_{t}w(t;z)}{w(t;z)-w(x;z)}-
          \frac{\partial_{t}w(t;z)}{w(t;z)-w(0;z)}=\\
          \frac{x}{t(t-x)}+\frac{1}{t-X_{1}(Y_{0}(x;z);z)}+\frac{1}{t-X_{1}(Y_{1}(x;z);z)}-\frac{1}{t+1},
     \end{multline*}
we immediately obtain the following.
     \begin{cor}\label{v-explicit_h(x;z)_second}
     The function $h^{i_{0},j_{0}}(x;z)$ can be split as $h^{i_{0},j_{0}}(x;z)=h_{1}^{i_{0},j_{0}}(x;z)+h_{2}^{i_{0},j_{0}}(x;z)+h_{3}^{i_{0},j_{0}}(x;z)$, where
         \begin{eqnarray}
              h_{1}(x;z)\hspace{-1mm}&=& x^{i_{0}}Y_{0}(x;z)^{j_{0}},\\
              h_{2}(x;z)\hspace{-1mm}&=& \frac{x}{\pi }
              \int_{x_{1}(z)}^{x_{2}(z)}
              \frac{t^{n_{0}-1}}{t-x}\mu_{j_{0}}
              (t;z)[-d(t;z)]^{1/2}\textnormal{d}t,
              \label{v-new_integral_form_h2}\\
              h_{3}(x;z)\hspace{-1mm}&=& \frac{1}{\pi }
              \int_{x_{1}(z)}^{x_{2}(z)}
              t^{i_{0}}\left[\frac{1}{t-X_{1}(Y_{0}(x;z);z)}+
              \right.\label{v-new_integral_form_h} \\
              &&\quad\quad\quad\quad\quad\hspace{2.3mm}\left.\frac{1}{t-X_{1}(Y_{1}(x;z);z)}
              -\frac{1}{t+1}\right]\mu_{j_{0}}(t;z)
              [-d(t;z)]^{1/2}\textnormal{d}t\nonumber.
          \end{eqnarray}
     \end{cor}

The end of Section \ref{DistributionTheVoterModel} aims at obtaining an
expression of $h^{i_{0},j_{0}}(1;z)$ which is efficient---in the
sense of computing the asymptotic of its coefficients (that we shall do in Section \ref{AsymptoticTheVoterModel}). In order
to achieve this, we shall make the change of variable $\widehat{b}(t;z)
=b(t;z)/[4a(t;z)c(t;z)]^{1/2}$ in the integrals \eqref{v-new_integral_form_h2} and \eqref{v-new_integral_form_h} of Corollary
\ref{v-explicit_h(x;z)_second}. The main reason of this is that using \eqref{v-def_mu} yields
     \begin{equation}\label{v-natural_expression_mu}
          \mu_{j_{0}}(t;z)[-d(t;z)]^{1/2}=
          \left(\frac{c(t;z)}{a(t;z)}\right)^{j_{0}/2}
          U_{j_{0}-1}(-\widehat {b}(t;z))
          [1-\widehat{b}(t;z)^{2}]^{1/2},
     \end{equation}
where the $(U_{n})_{n\in {\bf Z}_+}$ are the Chebyshev polynomials
of the second kind. {We refer to \cite{SZ} for a complete exposition on these polynomials.} We recall that they are the orthogonal
polynomials associated with the weight $t\mapsto
[1-t^{2}]^{1/2}{\bf 1}_{]-1,1[}(t)$, and that their
explicit expression is
\begin{equation*}
U_{n}(u)=
          \frac{(u+[u^{2}-1]^{1/2})^{n+1}-
          (u-[u^{2}-1]^{1/2})^{n+1}}{2[u^{2}-1]^{1/2}},
          \quad \forall u\in {\bf C},\quad \forall n\in{\bf Z}_+.\end{equation*}
We also recall two properties of the Chebyshev polynomials of the second
kind that we will especially use here (see \cite{SZ} for their proof):
\begin{itemize}
\item They have the parity
of their order, in other words, for all
$u\in{\bf C}$ and all $n\in{\bf Z}_+$, $U_{n}(-u)=
(-1)^{n}U_{n}(u)$;
\item Their expansion in
the neighborhood of $1$ is $U_{n}(u)=(n+1)[1+n(n+2)(u-1)/3
+O(u-1)^2]$.
\end{itemize}
Further, function $t\mapsto \widehat{b}(t;z)$ is clearly a diffeomorphism between $]x_{1}(z),x_{2}(z)[$ and
$]-1,1[$; in addition, $\widehat{b}(t;z)=u$ implies $b(t;z)^{2}-4u^{2}a(t;z)c(t;z)=0$, which,
as a polynomial in the variable $t$, is reciprocal, so that we can quite easily obtain
and write the explicit expression of its roots, called the $t_{\ell}(u;z)$, $\ell\in\{1,\ldots,4\}$.
Defining $T(u;z)=3/z+u^{2}-u[2+u^2+6/z]^{1/2}$, then $t_{2}(u;z)=T(u;z)-[T(u;z)^{2}-1]^{1/2}$,
$t_{3}(u;z)=T(u;z)+[T(u;z)^{2}-1]^{1/2}$, $t_{1}(u;z)=t_{2}(-u;z)$ and $t_{4}(u;z)=t_{3}(-u;z)$.
Notice that we have enumerated the $t_{\ell}(u;z)$ in such a way that $t_{\ell}(1;z)=x_{\ell}(z)$ for any $\ell\in\{1,\ldots ,4\}$.
Moreover, it turns out that for $u\in ]-1,1[$, $\widehat{b}(t_{2}(u;z);z)=-u$,
so that the following result is an immediate consequence of the change of
variable $t=t_{2}(u;z)$ in Corollary \ref{v-explicit_h(x;z)_second} as well
as of the identity \eqref{v-natural_expression_mu}.

\begin{cor}
\label{v-explicit_h(x;z)_third}
We have $h^{i_{0},j_{0}}(1;z)=h_{1}^{i_{0},j_{0}}(1;z)+
     h_{2}^{i_{0},j_{0}}(1;z)+h_{3}^{i_{0},j_{0}}(1;z)$, where
     \begin{eqnarray*}
          h_{1}^{i_{0},j_{0}}(1;z)&=&\left(\frac{(3-z-3[(1-z)
          (1+z/3)]^{1/2}}{(2z)}\right)^{j_{0}}, \\
          h_{2}^{i_{0},j_{0}}(1;z)&=&\frac{1}{\pi}\int_{-1}^{1}
          \frac{U_{j_{0}-1}(u) t_{2}(u;z)^{i_{0}+j_{0}/2-1}}
          {t_{2}(u;z)-1}\partial_{u} t_{2}(u;z)
          [1-u^{2}]^{1/2}\textnormal{d}u,\phantom{\widetilde{\widetilde{\widetilde{\widetilde{\widetilde{I}}}}}}\\
          h_{3}^{i_{0},j_{0}}(1;z)&=&\frac{1}{\pi }
          \int_{-1}^{1}U_{j_{0}-1}(u)t_{2}(u;z)^{i_{0}+j_{0}/2}
          \left[\frac{1}{t_{2}(u;z)-X_{1}(Y_{0}(1;z);z)}+\right.\\
          &&\hspace{8mm}\left.+
          \frac{1}{t_{2}(u;z)-X_{1}(Y_{1}(1;z);z)}-
          \frac{1}{t_{2}(u;z)+1}\right]
          \partial_{u} t_{2}(u;z)
          [1-u^{2}]^{1/2}\textnormal{d}u.
     \end{eqnarray*}
\end{cor}

\section{Asymptotic tail distribution of the hitting times}
\label{AsymptoticTheVoterModel}
\setcounter{equation}{0}

{In this section we prove Theorem \ref{main_result_group_six_drift_zero}, which is the main result in this paper.} Let $\mathscr{D}$ denote the open unit disc: $\mathscr{D}=\{z\in {\bf C}: |z|<1\}$. In order to prove Theorem \ref{main_result_group_six_drift_zero}, we are going to prove that $h^{i_{0},j_{0}}(1;z)$ is holomorphic in $(1+\epsilon)\mathscr{D}\setminus [1,1+\epsilon[$ and that in the neighborhood of $1$,
     \begin{equation}
     \label{v-before_to_conclude}
          h^{i_{0},j_{0}}(1;z)=(3/4)3^{1/2}i_{0}j_{0}(i_{0}+j_{0})
          [1-z]^{3/2}[1+o(1)]+h_{0}^{i_{0},j_{0}}(z),
     \end{equation}
where $h_{0}^{i_{0},j_{0}}$ is holomorphic at $1$; it will then be enough to use classical singularity analysis (see \cite[Sections 6.2--6.4]{FLAJ}).

For this, according to Corollary \ref{v-explicit_h(x;z)_third}, we shall
consider successively $h_{1}^{i_{0},j_{0}}(1;z)$, $h_{2}^{i_{0},j_{0}}(1;z)$ and
$h_{3}^{i_{0},j_{0}}(1;z)$ in Proposition \ref{v-lemma_h1}. Equation \eqref{v-before_to_conclude} and Theorem \ref{main_result_group_six_drift_zero}
will then be a direct consequence of these three results.

\begin{prop}\label{v-lemma_h1}
The functions $h_{1}^{i_{0},j_{0}}(1;z)$, $h_{2}^{i_{0},j_{0}}(1;z)$ and $h_{3}^{i_{0},j_{0}}(1;z)$ are holomorphic in $(1+\epsilon)\mathscr{D}\setminus [1,1+\epsilon[$. Moreover, in the neighborhood of $1$, we have
     \begin{align*}
          h_{1}^{i_{0},j_{0}}(1;z)&=-j_{0}3^{1/2}[1-z]^{1/2}[1+
          (3+4{j_{0}}^{2})(1-z)/8+f_{1,1}^{i_{0},j_{0}}(z)(z-1)^{2}]
          +f_{1,2}^{i_{0},j_{0}}(z),\\
           h_{2}^{i_{0},j_{0}}(1;z)&=\frac{3^{1/2}j_{0}}{2}[1-z]^{1/2}
          [1+(1/2)(3/4+{j_{0}}^{2})(1-z)
          +f_{2,1}^{i_{0},j_{0}}(z)(1-z)^{2}]\\&+
          \frac{3^{1/2}j_{0}}{2\pi }(i_{0}+j_{0}/2-
          1/2)\ln(1-z)[1+(1-z)
          f_{2,2}^{i_{0},j_{0}}(z)]+f_{2,3}^{i_{0},j_{0}}(z),\\
          h_{3}^{i_{0},j_{0}}(1;z)&=\frac{3^{1/2}j_{0}}{16}(1-z)^{1/2}
          [8+(3+4{j_{0}}^{2}+12i_{0}(i_{0}+j_{0}))(1-z)
          +f_{3,1}^{i_{0},j_{0}}(z)(1-z^{2})]\\&-
          \frac{3^{1/2}j_{0}}{4\pi }(2i_{0}+j_{0}-
          1)\ln(1-z)[1+(1-z)f_{3,2}^{i_{0},j_{0}}(z)]
          +f_{3,3}^{i_{0},j_{0}}(z),
     \end{align*}
where the $f_{k,\ell}^{i_{0},j_{0}}$ are holomorphic at $1$.
\end{prop}

\begin{proof}
The proof of the facts dealing with $h_{1}^{i_{0},j_{0}}(1;z)$ directly
follows from the expression of this function written in
Corollary \ref{v-explicit_h(x;z)_third}.

Let us now focus on $h_{2}^{i_{0},j_{0}}(1;z)$. We recall from Corollary \ref{v-explicit_h(x;z)_third} that
     \begin{equation}
     \label{v-explicit_expression_h2}
          h_{2}^{i_{0},j_{0}}(1;z)=\frac{1}{\pi}\int_{-1}^{1}
          \frac{U_{j_{0}-1}(u) t_{2}(u;z)^{i_{0}+j_{0}/2-1}}
          {t_{2}(u;z)-1}\partial_{u} t_{2}(u;z)
          [1-u^{2}]^{1/2}\textnormal{d}u,     \end{equation}
where $t_{2}(u;z)=T(u;z)-[T(u;z)^{2}-1]^{1/2}$ and $T(u;z)=3/z+u^{2}-u[2+u^{2}+6/z]^{1/2}$.
In particular, the fact that $h_{2}^{i_{0},j_{0}}(1;z)$ is holomorphic in $(1+\epsilon)
\mathscr{D}\setminus [1,1+\epsilon[$ is clear, since making the change of variable
$u\mapsto -u$ in \eqref{v-explicit_expression_h2} allows us to write it as the
integral on $[0,1]$ of some function of $(u,z)$ holomorphic in $\mathscr{D}\times ((1+\epsilon)\mathscr{D}
\setminus [1,1+\epsilon[)$---note that although $T(u;z)$ has algebraic singularities, any function symmetrical in
$(T(u;z),T(-u;z))$ is meromorphic w.r.t.\ the variable $z$.

We now study the behavior of $h_{2}^{i_{0},j_{0}}(1;z)$ in the
neighborhood of $1$.
For this, we first transform \eqref{v-explicit_expression_h2},
until obtaining an expression that makes clearly appear the singularities
of $h_{2}^{i_{0},j_{0}}(1;z)$.

An easy calculation entails that $\partial_{u} t_{2}=\partial_{u}T/(1-t_{3}^2)$.
Moreover, by definition of the $t_{\ell}$ (see Section \ref{DistributionTheVoterModel}),
$(z^{2}/36)\prod_{\ell=1}^{4}(t-t_{\ell}(u;z))$ is equal to $b(t;z)^{2}-4u^{2}a(t;z)c(t;z)$.
In particular, $\prod_{\ell=1}^{4}(1-t_{\ell}(u;z))=(36/z^{2})(1-z(1+2u)/3)(1-z(1-2u)/3)$.
So we have
     \begin{equation}\label{v-simplif}
          \frac{\partial_{u}t_{2}(u;z)
          }
          {t_{2}(u;z)-1}=\frac{z^{2}
          \partial_{u}T(u;z)(1-t_{1}(u;z))
          (1-t_{4}(u;z))
          (1-t_{2}(u;z))
          }{18(1-z(1-2u)/3)
          (t_{2}(u;z)-t_{3}(u;z))
          (1-z(1+2u)/3)}.
     \end{equation}

We now expand the quantity $(1-t_{2}(u;z))t_{2}(u;z)^{i_{0}+j_{0}/2-1}$
according to the powers of $[T(u;z)^{2}-1]^{1/2}$, say $(1-t_{2}(u;z))
t_{2}(u;z)^{i_{0}+j_{0}/2-1}=\sum_{k\geq 0}F_{k}^{i_{0},j_{0}}(u;z)
[T(u;z)^{2}-1]^{k/2}$.
With these notations, \eqref{v-explicit_expression_h2} and
\eqref{v-simplif}, we obtain
     \begin{multline}
     \label{v-h2_as_sum_integrals}
          h_{2}^{i_{0},j_{0}}(1;z)=\sum_{k\geq 0}\int_{-1}^{1}\frac{z^{2}
          \partial_{u}T(u;z)(1-t_{1}(u;z))
          (1-t_{4}(u;z))}{18(1-z(1-2u)/3)}
          F_{k}^{i_{0},j_{0}}(u;z)\times \\
          \times \frac{[T(u;z)^{2}-1]^{k/2}}
          {(t_{2}(u;z)-t_{3}(u;z))
          (1-z(1+2u)/3)}
          U_{j_{0}-1}(u)[1-u^{2}]^{1/2}\text{d}u.
     \end{multline}
{In what follows, we analyze the behavior at $1$ of each integral in the sum \eqref{v-h2_as_sum_integrals}, first for $k\in\{0,1,2\}$, then for $k\geq 3$.}

\subsection*{Integrals corresponding to $k\in\{0,1,2\}$ in the sum \eqref{v-h2_as_sum_integrals}}
First, note that
     \begin{eqnarray*}
          F_{0}^{i_{0},j_{0}}&=&T^{i_{0}+j_{0}/2-1}
          (1-T),\\
           F_{1}^{i_{0},j_{0}}&=&T^{i_{0}+j_{0}/2-2}
          [T-(i_{0}+j_{0}/2-1)(1-T)],\\
          F_{2}^{i_{0},j_{0}}&=&T^{i_{0}+j_{0}/2-3}
          (i_{0}+j_{0}/2-1)[(1-T)(i_{0}+j_{0}/2-2)/2-T].
     \end{eqnarray*}

Now we set $F^{j_{0}}(u;z)=-z^{2}\partial_{u}T(u;z)(1-t_{1}(u;z))(1-t_{4}(u;z))
U_{j_{0}-1}(u)/(36(1-z(1-2u)/3))$, as well as,
     \begin{eqnarray*}
          G_{0}^{i_{0},j_{0}}(u;z)\hspace{-1.5mm}&=&\hspace{-1.5mm} [F^{j_{0}}(u;z)
          F_{0}^{i_{0},j_{0}}(u;z)z^{2}[T(-u;z)^{2}-1]^{1/2}]/
          [3(z+3)(1-z(1-2u)/3)^{1/2}],\\
          G_{1}^{i_{0},j_{0}}(u;z)\hspace{-1.5mm}&=&\hspace{-1.5mm} \phantom{[}
          F^{j_{0}}(u;z)F_{1}^{i_{0},j_{0}}(u;z),\\
          G_{2}^{i_{0},j_{0}}(u;z)\hspace{-1.5mm}&=&\hspace{-1.5mm} [F^{j_{0}}(u;z)
          F_{2}^{i_{0},j_{0}}(u;z)3(z+3)
          (1-z(1-2u)/3)^{1/2}]/[
          z^{2}[T_{1}(-u;z)^{2}-1]^{1/2}].
     \end{eqnarray*}

Since $t_{2}(u;z)-t_{3}(u;z)=-2[T(u;z)^{2}-1]^{1/2}$ and since
$(t_{2}(u;z)-t_{3}(u;z))(t_{1}(u;z)-t_{4}(u;z))=12(z+3)^{2}[(1
-z(1+2u)/3)(1-z(1-2u)/3)]^{1/2}/z^{2}$, the sum of the three terms
for $k\in\{0,1,2\}$ in \eqref{v-h2_as_sum_integrals} is equal to
     \begin{equation}
     \label{starde}
          \sum_{k=0}^{2}\int_{-1}^{1}\frac{G_{k}^{i_{0},j_{0}}(u;z)
          [1-u^{2}]^{1/2}}{[1-z(1+2u)/3]^{(3-k)/2}}\text{d}u.
     \end{equation}

Using now the expansion of the Chebyshev polynomials at $1$
(see \cite{SZ}), we obtain the expansion
$G_{0}^{i_{0},j_{0}}(u;z)=-2j_{0}(u-1)/9-j_{0}(z-1)/3+\sum_{k+\ell\geq 2}
G_{0,k,\ell}^{i_{0},j_{0}}(u-1)^{k}(z-1)^{\ell}$.
Then, with a repeated use of \eqref{v-second_elliptic_integral}
of Lemma \ref{v-asymptotic_integrals_lemma_2},
we get
     \begin{equation*}
          \int_{-1}^{1}\frac{G_{0}^{i_{0},j_{0}}(u;z)[1-u^{2}]^{1/2}}
          {[1-z(1+2u)/3]^{3/2}}\text{d}u=
          j_{0}3^{1/2}\ln(1-z)[(1-z)/4+(1-z)^{2}g_{0}^{i_{0},j_{0}}(z)]+f_{0}^{i_{0},j_{0}}(z),
     \end{equation*}
$f_{0}^{i_{0},j_{0}}$ and $g_{0}^{i_{0},j_{0}}$ being holomorphic at $1$.

In the same way, $G_{1}^{i_{0},j_{0}}(u;z)=\sum_{k+\ell\geq 2}G_{1,k,\ell}^{i_{0},j_{0}}
(u-1)^{k}(z-1)^{\ell}-j_{0}/3-j_{0}(6{j_{0}}^{2}+35-48i_{0}-$
$24j_{0})(u-1)/54
+j_{0}(-53+48i_{0}+24j_{0})(z-1)/36$. A repeated application of Lemma
\ref{v-asymptotic_integrals_lemma_1} then gives that
     \begin{eqnarray*}
          \int_{-1}^{1}\frac{G_{1}^{i_{0},j_{0}}(u;z)[1-u^{2}]^{1/2}}
          {1-z(1+2u)/3}\text{d}u\hspace{-2mm}&=&\hspace{-2mm}
          f_{1}^{i_{0},j_{0}}(z)+j_{0}3^{1/2}[1-z]^{1/2}\times
          \\\hspace{-2mm}&\times &\hspace{-2mm}
          [1/2+(3/4+{j_{0}}^{2})(1-z)/4+(1-z)^{2}g_{1}(z)],
     \end{eqnarray*}
$f_{1}^{i_{0},j_{0}}$ and $g_{1}^{i_{0},j_{0}}$ being holomorphic at $1$.

Finally, we have $G_{2}^{i_{0},j_{0}}(u;z)=2j_{0}(i_{0}+j_{0}/2-1)/3+
\sum_{k+\ell\geq 1}G_{2,k,\ell}^{i_{0},j_{0}}(u-1)^{k}(z-1)^{\ell}$.
So with a repeated use of \eqref{v-first_elliptic_integral}
of Lemma \ref{v-asymptotic_integrals_lemma_2}, we reach the conclusion that
     \begin{eqnarray*}
          \int_{-1}^{1}\frac{G_{2}^{i_{0},j_{0}}(u;z)[1-u^{2}]^{1/2}}
          {[1-z(1+2u)/3]^{1/2}}\text{d}u\hspace{-2mm}&=&\hspace{-2mm}
          f_{2}^{i_{0},j_{0}}(z)+j_{0}(i_{0}+j_{0}/2-1)3^{1/2}\times
          \\\hspace{-2mm}&\times&\hspace{-2mm}
          \ln(1-z)[(1-z)/2+(1-z)^{2}g_{2}^{i_{0},j_{0}}(z)],
     \end{eqnarray*}
$f_{2}^{i_{0},j_{0}}$ and $g_{2}^{i_{0},j_{0}}$ being holomorphic at $1$.

\subsection*{Integrals corresponding to $k\geq 3$ in the sum \eqref{v-h2_as_sum_integrals}}

Note first that if $k$ is odd and larger than $3$, the
associated function in \eqref{v-h2_as_sum_integrals} is in fact
holomorphic in the neighborhood of $1$: indeed, for this it is
enough to notice that $t_{2}(u;z)-t_{3}(u;z)=-2[T(u;z)^{2}-1]^{1/2}$.
For this reason, all the terms in \eqref{v-h2_as_sum_integrals} associated
with values of $k$ which are odd and larger than $3$ do not have any singularity at $1$.

On the other hand, if $k$ is even and larger than $3$, then the underlying
term in the sum \eqref{v-h2_as_sum_integrals} can be written as
     \begin{equation*}
          \int_{-1}^{1}[1-z(1+2u)/3]^{(k-3)/2}H_{k}^{i_{0},j_{0}}(u;z)[1-u^{2}]^{1/2}\text{d}u,
     \end{equation*}
where the function $H_{k}^{i_{0},j_{0}}(u;z)$ is holomorphic in the
neighborhood of $(1,1)$. The last integral is obviously equal to
     \begin{equation*}
          \int_{-1}^{1}[1-z(1+2u)/3]^{(k-2)/2}H_{k}^{i_{0},j_{0}}
          (u;z)[1-u^{2}]^{1/2}[1-z(1+2u)/3]^{-1/2}\text{d}u.
     \end{equation*}

Then, expanding $[1-z(1+2u)/3]^{(k-2)/2}H_{k}^{i_{0},j_{0}}(u;z)$ w.r.t.\ the
powers of $(u-1)^{k}(z-1)^{\ell}$ and using \eqref{v-first_elliptic_integral} of
Lemma \ref{v-asymptotic_integrals_lemma_2}, we obtain that the integral above
equals $\ln(1-z)(z-1)^{k-2}g_{k}^{i_{0},j_{0}}(z)+f_{k}^{i_{0},j_{0}}(z)$,
$f_{k}^{i_{0},j_{0}}$ and $g_{k}^{i_{0},j_{0}}$ being holomorphic at $1$.

Finally, the sum of all the terms corresponding in \eqref{v-h2_as_sum_integrals} to values of $k$ which are even and
larger than $3$ can be written, in the neighborhood of $1$, as $\ln(1-z)(1-z)^{2}g^{i_{0},j_{0}}(z)+f^{i_{0},j_{0}}(z)$,
where $f^{i_{0},j_{0}}$ and $g^{i_{0},j_{0}}$ are holomorphic at $1$.

\subsection*{End of the proof of Proposition \ref{v-lemma_h1}}

Putting the latter fact together with \eqref{starde}
 concludes the proof of the expansion for $ h_{2}^{i_{0},j_{0}}(1;z)$ stated in Proposition \ref{v-lemma_h1}. Finally, the proof of the facts regarding $h_{3}^{i_{0},j_{0}}(1;z)$,
via a repeated use of Lemmas \ref{v-asymptotic_integrals_lemma_1}
and \ref{v-asymptotic_integrals_lemma_2}, is totally similar to that for $h_{2}^{i_{0},j_{0}}(1;z)$,
 so we omit it.
\end{proof}

\appendix

\section{}
{In this appendix, we state and prove two technical lemmas, which concern the behavior of some integrals with parameters near their singularities. These lemmas are crucial for the proof of our main results, but are independent of the rest of the paper.}

\begin{lem}
\label{v-asymptotic_integrals_lemma_1}
For any $k\in{\bf Z}_{+}$, let $P_{k}$ be the principal part at infinity of
$[Z^{2}-1]^{1/2}(1-Z)^{k}$, i.e., the unique polynomial such that $[Z^{2}
-1]^{1/2}(1-Z)^{k}-P_{k}(Z)$ goes to zero as $|Z|$ goes to infinity. Then
     \begin{equation*}
          \int_{-1}^{1}\frac{(1-u)^{k}[1-u^{2}]^{1/2}}{1-z(1+2u)/3}
          \textnormal{d}u =\frac{3\pi}{2z}\left[(1+z/3)^{1/2}\left(\frac{-3}{2z}\right)^{k+1}
          (1-z)^{k+1/2}+P_{k}\left(\frac{3}{2z}-\frac{1}{2}\right)\right].
     \end{equation*}
\end{lem}

\begin{proof}[Proof]
For $\epsilon>0$, we consider the closed contour $\mathscr{A}_{\epsilon}^{+}\cup \mathscr{A}_{\epsilon}^{-}\cup
\mathscr{B}_{\epsilon}^{+}\cup \mathscr{B}_{\epsilon}^{-}$, where $\mathscr{A}_{\epsilon}^{\pm}=\{\pm 1\mp
i \epsilon \exp(i t),t\in [0,\pi]\}$ and $\mathscr{B}_{\epsilon}^{\pm}=\{\pm i \epsilon \mp
t,t\in [-1,1]\}$. Then we apply on it the residue theorem at infinity to the function $(1-u)^{k}[1-u^{2}]^{1/2}
/[1-z(1+2u)/3]$ and we let $\epsilon $ going to zero.
\end{proof}

\begin{lem}\label{v-asymptotic_integrals_lemma_2}
For any $k\in{\bf Z}_{+}$, the integrals written in the left hand side
of \eqref{v-first_elliptic_integral} and \eqref{v-second_elliptic_integral} below
are holomorphic in $(1+\epsilon)\mathscr{D}\setminus[1,1+\epsilon[$ for
$\epsilon>0$ small enough. In the neighborhood of $1$, they are equal to
     \begin{eqnarray}
          \int_{-1}^{1}\frac{(1-u)^{k}[1-u^{2}]^{1/2}}
          {[1-z(1+2u)/3]^{1/2}}\textnormal{d}u
          &=&\ln(1-z)(1-z)^{k+1}\alpha_{k}(z)+\beta_{k}(z),
          \label{v-first_elliptic_integral}\\
          \int_{-1}^{1}\frac{(1-u)^{k}[1-u^{2}]^{1/2}}
          {[1-z(1+2u)/3]^{3/2}} \textnormal{d} u
          &=&\ln(1-z)(1-z)^{k}\gamma_{k}(z)+\delta_{k}(z),
          \label{v-second_elliptic_integral}
     \end{eqnarray}
where $\alpha_{k}$, $\beta_{k}$, $\gamma_{k}$ and $\delta_{k}$ are holomorphic at $1$,
$\alpha_{k}(1)\neq 0$ and $\gamma_{k}(1)\neq 0$. Furthermore, $\alpha_{0}(1)=3^{3/2}/4$,
$\gamma_{0}(1)=-3^{3/2}/2$, $\gamma_{0}'(1)=-3^{1/2}99/32$ and $\gamma_{1}(1)=3^{1/2}
27/8$.
\end{lem}

\begin{proof}[Proof]
The fact that the integrals written in the left hand side of \eqref{v-first_elliptic_integral}
and \eqref{v-second_elliptic_integral} are, for $\epsilon>0$ small enough, holomorphic in
$(1+\epsilon)\mathscr{D}\setminus[1,1+\epsilon[$ is clear from their expression.

Let us now study their behavior near $1$ and start by considering \eqref{v-first_elliptic_integral}.
Replace first the lower bound $-1$ by $-1/2$ in the integral \eqref{v-first_elliptic_integral}.
This is equivalent to add a function holomorphic in some $(1+\epsilon)\mathscr{D}$ and this will
eventually change $\beta_{k}$ but not $\alpha_{k}$ in the right hand side member of
\eqref{v-first_elliptic_integral}. Then, the change of variable $v^{2}=(1+2u)/3$ gives
     \begin{equation}\label{v-change_of_variable}
          \int_{-1/2}^{1}\frac{(1-u)^{k}[1-u^{2}]^{1/2}}
          {[1-z(1+2u)/3]^{1/2}}\textnormal{d}u
          =3^{1/2}\left(\frac{3}{2}\right)^{k+1}
          \int_{0}^{1}\frac{[1-v^{2}]^{k+1/2}}
          {[1-z v^{2}]^{1/2}}[1+3v^{2}]^{1/2} v \text{d}v.
     \end{equation}

By using the expansion of ${v}^{1/2}$ in the neighborhood of $1$,
we can develop the function $[1+3v^{2}]^{1/2}v$ according to the
powers of $v^{2}-1$: $[1+3v^{2}]^{1/2}v=2+(7/4)
[v^{2}-1]+\cdots $. Further, in \cite[Chapter F]{thesis},
we have proved, using the elliptic integrals theory, that for
any $k\in{\bf Z}_{+}$ there exist two functions $\phi_{k}$ and
$\psi_{k}$, holomorphic in the neighborhood of $1$ and satisfying
$\phi_{k}(1)\neq 0$, such that
     \begin{equation}
     \label{IA}
          \int_{0}^{1}\frac{[1-v^{2}]^{1/2+k}}
          {[1-z v^{2}]^{1/2}}\text{d}v
          =\ln(1-z)(z-1)^{k+1}\phi_{k}(z)+\psi_{k}(z),
     \end{equation}
we have there also proved that $\phi_{0}(1)=1/4$.
The equality \eqref{v-first_elliptic_integral} is then
an obvious outcome of \eqref{v-change_of_variable},
of the expansion of $[1+3v^{2}]^{1/2}v$ according to the
powers of $v^{2}-1$, and of the repeated use of \eqref{IA}. The fact that $\alpha_{0}(1)=3^{3/2}/4$ comes from
the equality $\phi_{0}(1)=1/4$.

Likewise, we prove the equality \eqref{v-second_elliptic_integral} and we obtain the announced values of $\gamma_{0}(1)$, $\gamma_{0}'(1)$ and $\gamma_{1}(1)$.
\end{proof}

\end{document}